\documentclass[11pt]{article}
\usepackage{amsmath, amsfonts, amssymb}
\usepackage{graphicx, amsthm}
\usepackage{bm}
\usepackage{xcolor}

\providecommand{\U}[1]{\protect\rule{.1in}{.1in}}

\newtheorem{theorem}{Theorem}
\newtheorem*{theoremA*}{Theorem A}
\newtheorem*{theoremB*}{Theorem B}
\newtheorem*{ass*}{Assumptions}

\newtheorem{corollary}[theorem]{Corollary}
\newtheorem{lemma}[theorem]{Lemma}
\newtheorem{proposition}[theorem]{Proposition}
\theoremstyle{definition}
\newtheorem{definition}[theorem]{Definition}
\newtheorem{remark}[theorem]{Remark}
\newtheorem{example}[theorem]{Example}

\begin{document}

\begin{center}
{\LARGE Desingularization of the Sweeping Process Mapping}
\end{center}

\medskip

\begin{center}
{\large \textsc{A. Daniilidis and S. Tapia-Garcia}}

\bigskip
\end{center}


\noindent\textbf{Abstract.} In \cite{DD}, the celebrated K\L -inequality has
been extended from definable functions $f:\mathbb{R}^{n}\rightarrow\mathbb{R}
$ to definable multivalued maps $S:\mathbb{R}\rightrightarrows\mathbb{R}^{n}$,
by establishing that the co-derivative mapping $D^{\ast}S$ admits a
desingularization around every critical value. As was the case in the gradient
dynamics, this desingularization yields a uniform control of the lengths of
all bounded orbits of the corresponding sweeping process $-\dot{\gamma}(t)\in
N_{S(t)}(\gamma(t))$. In this paper, working outside the framework of
o-minimal geometry, we characterize the existence of a desingularization for
the coderivative in terms of the behaviour of the sweeping process orbits and
the integrability of the talweg function. These results are close in spirit
with the ones in~\cite{BDLM}, where characterizations for the desingularization of the 
(sub)gradient of functions had been obtained.

\vspace{0.55cm}


\vspace{0.55cm}

\textit{Secondary} 26D10, 34A60, 37C10

\tableofcontents

\section{Introduction}

It is well-known that every $\mathcal{C}^{1}$ smooth function $f:\mathbb{R}
^{n}\rightarrow\mathbb{R}$ which is definable in some o-minimal structure has
finitely many critical values. 
Kurdyka \cite{Kurdyka} showed that if $\bar
{r}\in f(\mathbb{R}^{n})$ is a critical value and $\mathcal{U}$ is a nonempty
open bounded subset of $\mathbb{R}^{n},$ then there exist $\rho>0$ and a
$\mathcal{C}^{1}$-smooth function $\psi:[\bar{r},\bar{r}+\rho]\rightarrow
\lbrack0,+\infty)$ satisfying
\begin{equation}
\Vert\nabla(\psi\circ f)(x)\Vert\geq1,\qquad\text{for all }\mathcal{\;}
x\in\mathcal{U\;}\text{such that \ }f(x)\in(\bar{r},\bar{r}+\rho). \label{KL}%
\end{equation}
The above inequality generalizes to o-minimal functions the \L ojasiewicz
gradient inequality (established in \cite{L} for the class of $\mathcal{C}
^{1}$ subanalytic functions) and is nowadays known as the
Kurdyka-\L ojasiewicz inequality (in short, K\L -inequality). For definitions
and properties of o-minimal functions the reader is referred to \cite{vdD}.
Both the \L ojasiewicz and the K\L -inequality have been further extended to
nonsmooth (subanalytic and respectively o-minimal) functions, see \cite{BDL,
BDLS}. These inequalities allow to control uniformly the lengths of the
bounded\ (sub)gradient orbits, see \cite{Lojasiewicz1984,Kurdyka,BDL}
.\smallskip\newline One of the main features of Kurdyka's work \cite{Kurdyka}
was to consider the so-called \textit{talweg} function
\begin{equation}
m(r)=\,\underset{x\in\mathcal{U}}{\inf}\,\big\{\,\|\nabla f(x)\|:\,f(x)=r\big\},\qquad
r\in(\bar{r},\bar{r}+\rho), \label{talweg}
\end{equation}
which captures the worst behaviour (lower value of the norm of the gradient) at
the level set $[f=r].$ Kurdyka used the above function to defined the talweg
set $\mathcal{V}(r)$ consisting of points $x\in f^{-1}(r)$ with $||\nabla
f(x)||\leq\,2\,m(r).$ He then made use of a definable version of the \textit{curve selection
lemma} to obtain a smooth curve $r\mapsto\theta(r)\in\mathcal{V}(r) $ which is
directly linked to the \textit{desingularizing} function~$\psi.$ A
straightforward consequence of (\ref{KL}) is that the length of every bounded
gradient curve $\dot{\gamma}=-\nabla f(\gamma)$ contained in $f^{-1}((\bar
{r},\bar{r}+\rho))$ is majorized by $\psi(\bar{r}+\rho)-\psi(\bar{r})$ (and
therefore it is bounded). The same is true for the lengths of the piecewise
gradient curves, that is, curves obtained by concatenating countably many
gradient curves $\{\gamma_{i}\}_{i\geq1},$ where $\gamma_{i}\subset
f^{-1}([r_{i+1},r_{i}))$ and $\{r_{i}\}_{i}$ is a strictly decreasing sequence
in $(\bar{r},\bar{r}+\rho)$ converging to $\bar{r}$. These curves may have
countably many discontinuities. \smallskip\newline Outside the framework of
o-minimality the K\L -inequality (\ref{KL}) may fail even for $\mathcal{C}
^{2}$-smooth functions \cite[Section~4.3]{BDLM} or for $\mathcal{C}^{\infty}
$-smooth function with a unique critical value \cite[p. 12]{PdM1982}. Bolte,
Daniilidis, Ley and Mazet in \cite{BDLM} considered the problem of
characterizing the existence of a desingularization function~$\psi$ and the
validity of (\ref{KL}) for an upper isolated critical value $\bar{r}$ of a
semiconvex \textit{coercive} function $f$ defined in a Hilbert space. (A
function $f$ is called coercive, if it has bounded sublevel sets. This
assumption replaces the use of an open bounded set $\mathcal{U}$ in Kurdyka's
result.) We reproduce below one of the main results of the aforementioned
work, see \cite[Theorem 20]{BDLM}, for the special case where the function is
smooth and defined in finite dimensions.

\begin{theorem}
[characterization of the K\L -inequality]\label{theorem BDLM} Let
$f:\mathbb{R}^{n}\rightarrow\mathbb{R}\cup\{+\infty\}$ be a $\mathcal{C}^{2}
$-smooth (or more generally $\mathcal{C}^{1}$-smooth semi-convex) coercive
function and $\bar{r}\in f(\mathbb{R}^{n})$ an upper isolated critical value.
The following statements are equivalent:

\begin{enumerate}
\item[a)] \textbf{(K\L -inequality)} There exist $\rho>0$ and a smooth
function $\psi:[\bar{r},\bar{r}+\rho)\rightarrow\lbrack0,\infty)$ such that
\[
\Vert\nabla(\psi\circ f)(x)\Vert\geq1,~\text{for all }x\in f^{-1}((\bar
{r},\bar{r}+\rho)).
\]

\item[b)] \textbf{(Length control for gradient curves)} There exist $\rho>0$
and a strictly increasing continuous function $\sigma:[\bar{r},\bar{r}
+\rho)\rightarrow\lbrack0,\infty)$ with $\sigma(\bar{r})=0$ such that
\[
\int_{0}^{T}\Vert\dot{\gamma}(t)\Vert dt\leq\sigma(f(\gamma(0)))-\lim
_{t\rightarrow T}\sigma(f(\gamma(t))),\qquad\text{(we may have }
T=+\infty\text{)}
\]
for all gradient curves $\gamma:[0,T)\rightarrow\mathbb{R}^{n}$ with
$\gamma([0,T))\subset f^{-1}((\bar{r},\bar{r}+\rho)).$

\item[c)] \textbf{(Length bound for piecewise gradient curves)} There exist
$\rho,M>0$ such that
\[
\int_{0}^{T}\Vert\dot{\gamma}(t)\Vert dt\leq M,
\]
for all piecewise gradient curves $\gamma:[0,T)\rightarrow\mathbb{R}^{n}$ with
$\gamma([0,T))\subset f^{-1}((\bar{r},\bar{r}+\rho))$.

\item[d)] \textbf{(Integrability condition)} There exists $\rho>0$ such that
the function
\begin{equation*}
r\mapsto \,\underset{x\in f^{-1}(r)}{\sup }\dfrac{1}{\Vert \nabla f(x)\Vert }
,\qquad r\in (\bar{r},\bar{r}+\rho ),
\end{equation*}
is finite-valued and belongs to $\mathcal{L}^{1}(\overline{r},\overline{r}+\rho)$.
\end{enumerate}
\end{theorem}

\noindent Recently, Daniilidis and Drusvyatskiy \cite{DD} showed that every multivalued
map $S:\mathbb{R}\rightrightarrows\mathbb{R}^{n}$ with definable graph admits
a desingularization of its graphical coderivative $D^{\ast}S:\mathbb{R}
^{n}\rightrightarrows\mathbb{R}$ around any critical value $t\in\mathbb{R}$.
(Relevant definitions and a more precise statement are given in
Section~\ref{ss:DD}.) This result yields a uniform bound for the lengths of
all bounded orbits of the sweeping process defined by $S$ (see forthcoming
Definition~\ref{sweepdefn}). The aforementioned results of \cite{DD} are also
covering the results of Kurdyka in \cite{Kurdyka} by considering a sweeping
process mapping $S$ related to the sublevel sets of the smooth definable
function $f$ (\emph{c.f.} Remark~\ref{rem-1}).\smallskip\newline The main
contributions of this work are the following:

\begin{itemize}
\item Without assuming o-minimality, we characterize the desingularization of
the coderivative of a smooth sweeping process (see
Definition~\ref{smooth sweeping}) by establishing an analogous result to
Theorem~\ref{theorem BDLM}. This is the main result of this work, which is
resumed in Section~\ref{ss:thmA}.

\item Since the evolution of the sweeping process is not reversible in time,
we introduce in Definition~\ref{coderivative} an asymmetric version of the
modulus for the coderivative of a multivalued map $S$, $\Vert D^{\ast
}S(t,x)|^{+}$, that captures the orientation of the dynamics. (In \cite{DD},
the prevailing assumption of o-minimality made it possible to work directly
with the usual modulus.)

\item We establish an asymmetric version of \cite[Theorem 9.40]{RW} (sometimes
known as the Mordukhovich Criterion) relating the asymmetric modulus of the
coderivative to the \emph{oriented calmenss} of the multivalued map
(Proposition~\ref{oriented criterion}). We then obtain Theorem~B (Section~\ref{ss:thmB})
which relates the desingularization of the coderivative with the length of
discrete sequences given by the catching--up algorithm. (This algorithm can be
perceived as the proximal algorithm over a function $f$ whenever the
multivalued map $S$ is defined by the sublevel sets of $f$.)
\end{itemize}

\noindent The outline of this manuscript is as follows: In
Section~\ref{notation and preliminaries}, we fix our notation, we quote
preliminary results of variational analysis required in the sequel. In
Section~\ref{main results}, we fix our setting, explain our assumptions and
state the two main results of this paper (Theorem~A and Theorem~B). The proofs
of these results together with other auxiliary results will be given in
Section~\ref{Sec:Proofs}.

\section{Notation and Preliminaries}

\label{notation and preliminaries}The notation used along this paper is
standard and follows the lines of \cite{RW}. For any two nonempty sets
$A,~B\subset\mathbb{R}^{n}$, the excess of $A$ over $B$ is given by
$\mathrm{ex\,}{(A,B):=\sup\{d(x,B):~x\in A\}}$ and their Hausdorff-Pompeiu
distance is defined by $\mathrm{dist\,}(A,B):=\max\,\{\mathrm{ex\,}
(A,B),\,\mathrm{ex\,}(B,A)\}.$ \smallskip\newline Let $C\subseteq
\mathbb{R}^{n}$ be a closed set and let $x\in\mathbb{R}^{n}$. The set of
projections of $x$ at $C$ is defined by
$\mathrm{Proj}_C(x):=\{y\in C:\Vert x-y\Vert=d(x,C)\}$, where $d(x,C):=\inf_{y\in C}\,d(x,y)$. The Fr\'{e}chet normal
cone to $C$ at $x\in C$, denoted by $\hat{N}_{C}(x)$, is the set of vectors
$v\in\mathbb{R}^{n}$ satisfying
\[
\limsup_{\substack{y\in C \\y\rightarrow x}}\dfrac{\langle v,y-x\rangle}{\Vert
y-x\Vert}\leq0.
\]
The limiting normal cone to $C$ at $x$, denoted by $N_{C}(x)$, consists of all
vectors $v\in\mathbb{R}^{n}$ such that there exists a sequence $(x_{i}
)_{i}\subset C$ and $v_{i}\in\hat{N}_{C}(x_{i})$ satisfying $x_{i}\rightarrow
x$ and $v_{i}\rightarrow v$.\smallskip\newline

\subsection{Sweeping process dynamics}

Let $S:\mathbb{R}\rightrightarrows\mathbb{R}^{n}$ be a multivalued map. The
effective domain of $S$, denoted by $\mathrm{dom}(S)$, is the set
$\{t\in\mathbb{R}:~S(t)\neq\emptyset\}$. We denote by $\mathcal{S}
=\textup{gph}(S) $ the graph of the multivalued map $S,$ that is,
\[
\mathcal{S}=\mathrm{gph\,}(S):=\{(t,x)\in\mathbb{R}^{n+1}:~x\in S(t)\}.
\]

\noindent Let us introduce the following dynamical system, known as sweeping
process, determined by the multivalued function $S$. The definition
implicitely implies that $\mathrm{dom}(S)$ has nonempty interior, and is often
an interval (possibly unbounded). In particular, in our seeting (\textit{c.f} Assumptions in Section~\ref{ss:31}) 
 $\mathrm{dom}(S)$ will always be an interval (possibly unbounded).

\begin{definition}
[sweeping process dynamics]\label{sweepdefn} Let $S:\mathbb{R}
\rightrightarrows\mathbb{R}^{n}$ be a multivalued map and $I\subset
\mathrm{dom}(S)$ be a nonempty interval of $\mathbb{R}$. We say that the
absolutely continuous curve $\gamma:I\rightarrow\mathbb{R}^{n}$ is a solution
(orbit) of the \textit{sweeping process} defined by $S$ if
\begin{equation}
\begin{cases}
-\dot{\gamma}(t)\in N_{S(t)}(\gamma(t)),~\forall_{a.e.}~t\in I,\\
~\hspace{0.15cm}\gamma(t)\in S(t),\qquad\quad\,\,\forall t\in I,
\end{cases}
\label{sweepingequation}
\end{equation}
where $N_{S(t)}(\gamma(t))$ stands for the normal cone of $S(t)$ at
$\gamma(t)$.
\end{definition}

Notice that \eqref{sweepingequation} can be formally satisfied by curves with
possible discontinuities (the set of discontinuities has then to be of measure
zero). For our purposes it is useful to consider the class of
\textit{piecewise absolutely continuous }curves, that is, curves
$\gamma:I\rightarrow\mathbb{R}^{n}$ whose possible discontinuities are
contained in a closed countable set $D$ and being absolutely continuous on each interval of
$I\diagdown D.$ This latter set is open, therefore it is a countable union of
disjoint intervals $J_{i}$, and $\gamma$ is required to be absolutely
continuous on each $J_{i}$.

\smallskip

\textit{Notation} ($\mathcal{AC}(S,I)$, $\mathcal{PAC}(S,I)$). We denote by
$\mathcal{AC}(S,I)$ (respectively $\mathcal{PAC}(S,I)$) the set of absolutely
continuous (respectively, piecewise absolutely continuous) orbits of the
sweeping process $S$ defined on the interval $I\subset\mathrm{dom}(S).$ The
length of a (piecewise) absolutely continuous curve $\gamma:I\rightarrow
\mathbb{R}^{n}$ is given by the formula
\[
\ell(\gamma):=\int_{I}\Vert\dot{\gamma}(t)\Vert dt.
\]

\subsection{Coderivative, (oriented) modulus and (oriented) talweg.}

Let $S:\mathbb{R}\rightrightarrows\mathbb{R}^{n}$ be a multivalued map with
closed values.

\begin{definition}
[Coderivative]\label{coderivative}The (limiting) coderivative of $S$ at
$(t,x)\in\mathcal{S}$ in $u\in\mathbb{R}^{n}$ is defined as follows:
\[
D^{\ast}S(t,x)(u):=\{a\in\mathbb{R}:~(a,-u)\in N_{\mathcal{S}}(t,x)\}.
\]

\end{definition}

\noindent Therefore $D^{\ast}S(t,x):\mathbb{R}^{n}\rightrightarrows\mathbb{R}$ is a
multivalued map and
\[
(u,a)\in\mathrm{gph\,}D^{\ast}S(t,x)\quad\text{ if and only if\quad}(a,-u)\in
N_{\mathcal{S}}(t,x).
\]
Since $\mathrm{gph\,}D^{\ast}S(t,x)$ is a cone, the map $D^{\ast}S(t,x)$ is
positively homogeneous and we can define its modulus via the formula:
\[
\Vert D^{\ast}S(t,x)\Vert^{+}:=\,\underset{\Vert u\Vert\leq1}{\sup}\big\{
|a|:\,a\in D^{\ast}S(t,x)(u)\big\}  .
\]
Although the above definition of a modulus is classical and relates nicely to
the Lipschitz continuity of $S$ (\emph{c.f.} \cite[Theorem 9.40]{RW}), the
symmetry of the absolute value of $\mathbb{R}$ (representing the time in our
dynamics) does not fit to the non-reversible dynamics of the sweeping process.
To remedy this, one needs to replace $|a|$ in the above formula by
$a^{+}:=\max\{0,a\}$ which eventually gives rise to the following definition.

\begin{definition}
[Asymmetric modulus of coderivative]\label{Def:asym_mod}For every
$(t,x)\in\mathcal{S}$ we define the asymmetric modulus of the coderivative
$D^{\ast}S(t,x)$ as follows:
\[
\Vert D^{\ast}S(t,x)|^{+}=\sup\big\{a^{+}:~a\in D^{\ast}S(t,x)(u),~\Vert
u\Vert\leq1\big\},
\]
where we adopt the convention $\sup(\emptyset)=0$.
\end{definition}

The following example gives some insight about the difference between the two moduli.

\begin{example}
\label{example 1} Let $f:\mathbb{R}^{n}\rightarrow\mathbb{R}$ be a
$\mathcal{C}^{1}$-smooth function and set
\[
S(r)=[f\leq r]:=\{x\in\mathbb{R}^{n}:~f(x)\leq r\},\qquad\text{for all }
r\in\mathbb{R}\text{.}
\]
This defines a multivalued map $S:\mathbb{R}\rightrightarrows\mathbb{R}^{n}$
associated to $f$, in the sense that the graph $\mathcal{S}$ of $S$ is the epigraph of~$f$ (up to a permutation of coordinates that brings the first coordinate of $\mathbb{R}^{n+1}$ to the last position). Let
$x\in S(r)$.\smallskip\newline If $f(x)<r$, then $x\in\mathrm{int}(S(r))$ and
$N_{\mathcal{S}}(r,x)=\{0\},$ yielding $\Vert D^{\ast}S(r,x)\Vert^{+}=\Vert
D^{\ast}S(r,x)|^{+}=0$. 
On the other hand, since the normal space of
$\mathrm{gph}(f)$ at $(x,f(x))$ is exactly $\mathbb{R}(\nabla f(x),-1)$, if
$f(x)=r$, then $N_{\mathcal{S}}(r,x)=\mathbb{R}_{+}(-1,\nabla f(x))$. Thus,
\[
\Vert D^{\ast}S(t,x)\Vert^{+}=\dfrac{1}{\Vert\nabla f(x)\Vert},~\hspace
{0.5cm}\text{but }~\hspace{0.5cm}\Vert D^{\ast}S(t,x)|^{+}=0.
\]

\end{example}

We now define the \emph{oriented talweg} function associated to the
multivalued map $S:\mathbb{R}\rightrightarrows\mathbb{R}^{n}.$ This captures
the worst case (larger value of the oriented modulus of the coderivative) on
each set $S(t),$ $t\in\mathbb{R}.$ This function will play an important role
in our main result.

\begin{definition}
[oriented talweg]\label{rem-talweg}The oriented talweg function of $S$ denoted
by ${\varphi^{\uparrow}}$ is defined as follows:
\[
{\varphi^{\uparrow}}(t)=\,\underset{x\in S(t)}{\sup}\{\Vert D^{\ast
}S(t,x)|^{+}\},\quad\text{for all }t\in\mathrm{dom}(S).
\]

\end{definition}

\begin{remark}
[Asymmetric structures]\label{rem-ven} In \cite{DD} the usual talweg function
$\varphi$ has been considered, based on the (symmetric) modulus of the
coderivative.%
\[
\varphi(t)=\,\underset{x\in S(t)}{\sup}\{\Vert D^{\ast}S(t,x)\Vert^{+}
\},\quad\text{for all }t\in\mathrm{dom}(S).
\]
The difference between $\varphi$ and ${\varphi^{\uparrow}}$ is that the modula
$\Vert D^{\ast}S(t,x)\Vert^{+},$ $(t,x)\in\mathcal{S},$ are now replaced by
their asymmetric versions $\Vert D^{\ast}S(t,x)|^{+}.$ The reader might notice
that $a^{+}:=\max\{0,a\}$ is a typical asymmetric norm of $\mathbb{R}$ and
$\Vert D^{\ast}S(t,x)|^{+}$ can be seen as a natural asymmetrization of the
modulus $\Vert D^{\ast}S(t,x)\Vert^{+}.$ The use of asymmetric objects seems
to be a natural tool in nonsmooth dynamics as well as in operations research
(orientable graphs). More details on asymmetric structures can be found in
\cite{Cobzas} and \cite{DV}.
\end{remark}

\subsection{Desingularization of the coderivative (definable case).}

\label{ss:DD}We now recall the main result of \cite{DD}. If $S:\mathbb{R}
\rightrightarrows\mathbb{R}^{n}$ is a multivalued map with a closed \textit{bounded}
graph $\mathcal{S}$, then assuming that $\mathcal{S}$ is definable in some
o-minimal structure, for every $a\in\mathbb{R}$ there exists $\rho>0$ and a
strictly increasing, continuous function $\Psi\colon\lbrack0,\rho
]\rightarrow\mathbb{R}$ that is $\mathcal{C}^{1}$-smooth on $(0,\rho)$, it
satisfies $\Psi(0)=a$ and $\Psi^{\prime}(r)>0$ for all $r\in(0,\rho)$ and
\begin{equation}
\| D^{\ast}(S\circ\Psi)(r,x) \|^{+}\leq1\qquad\text{ for
all }r\in(0,\rho)\text{ and all }x\in S(\Psi(r)). \label{eq:psi}
\end{equation}
It is easily seen that $\Psi$ is a homeomorphism between $[0,\rho]$ and
$[a,b]$ where $b=\Psi(\rho)$ and a diffeomorphism between $(0,\rho)$ and
$(a,b).$ Inequality (\ref{eq:psi}) has a particular interest when
$a\in\mathbb{R}$ is a \emph{critical value} of the coderivative $D^{\ast}S$ of
the sweeping process, that is,
\[
\varphi(t)=\underset{x\in S(t)}{\sup}\,\|D^{\ast}S(t,x)\|^{+}=+\infty.
\]
In this case we say that $\Psi$ \emph{desingularizes} the (modulus of the
coderivative around the) critical value $a$. The assumption of o-minimality on $S$ 
guarantees that the set of critical values is finite. In \cite{DD} it has
further been established, as consequence of \eqref{eq:psi}, that all bounded
orbits of the sweeping process $S$ have finite length and that the talweg
function $\varphi$ is integrable on $[a,b].\smallskip$\newline 
Let us notice that $\Vert D^{\ast}S(t,x)|^{+}\leq\Vert D^{\ast
}S(t,x)||^{+}$ (and consequently ${\varphi^{\uparrow}}(t)\leq\varphi(t)$) for all $t\in\lbrack a,b)$ and $x\in S(t)$. 
Therefore, we obtain the following.
\begin{corollary}
[Desingularization of oriented coderivative -- definable case]
\label{cor:des-as}If $S:\mathbb{R}\rightrightarrows\mathbb{R}^{n}$ is a
multivalued map with a closed definable bounded graph, then for every $a\in\mathbb{R}$
(possibly critical for the oriented modulus) there exists $\rho>0$ and $b>a$
such that:\smallskip\newline\textrm{(i).} there exists an increasing
homeomorphism $\Psi\colon\lbrack0,\rho]\rightarrow\lbrack a,b]$ which is
$\mathcal{C}^{1}$-diffeomorphism on $(0,\rho)$ such that:
\begin{equation}
||D^{\ast}(S\circ\Psi)(r,x)|^{+}\leq1\qquad\text{ for all }r\in(0,\rho)\text{
and all }x\in S(\Psi(r)). \label{eq:psi-as}
\end{equation}
\newline\textrm{(ii).} $\int_{a}^{b} {\varphi^{\uparrow}}(t)<\infty$ (the
oriented talweg function is integrable).
\end{corollary}

\begin{remark}\label{remark 9}
[Relation with the K\L -inequality]\label{rem-1} (i). The described
desingularization of the coderivative can be seen as a generalization of the
K\L -inequality for $\mathcal{C}^{1}$-smooth definable functions (established
by Kurdyka in \cite{Kurdyka}) in the following sense: let $f:\mathbb{R}^{n}\rightarrow\mathbb{R}$ 
be a $\mathcal{C}^{1}$-smooth coercive function which is definable in some o-minimal structure. Then, the multivalued function
\begin{equation}
\left\{
\begin{array}
[c]{l}
S_{f}:\mathbb{R}\rightrightarrows\mathbb{R}^{n}\smallskip\\
S_{f}(t)=[f\leq-t],\quad t\in\mathbb{R}
\end{array}
\right.  \label{eq:Sf}
\end{equation}
is o-minimal (it is definable in the same o-minimal structure as $f$) and the
desingularization of its gradient described in (\ref{KL}) can be deduced from
the desingularization coderivative of $S$ and vice versa. We refer the reader
to \cite[Section~5.1]{DD} for more details. \smallskip\newline
(ii). In~\cite{DD}, the assumption that $\mathcal{S}$ is bounded has not been made, and  similarly to \eqref{talweg}, the supremum of the definition of $\varphi(t)$ had to be taken over $S(t)\cap\mathcal{U}$, where $\mathcal{U}\subset\mathbb{R}^n$ is a fixed open bounded set, which gives rise to a talweg function 
$\varphi_{\mathcal{S}}$ depending on $\mathcal{U}$. Even if in Section~\ref{main results} we deal with potentially unbounded sweeping processes, we do not need to make use of $\mathcal{U}$, thanks to the assumptions given in Section~\ref{ss:31}. 
\end{remark}

\section{Characterization of desingularization of the coderivative}

\label{main results}

In this paper we are interested in sweeping process mappings $S$ that are not
o-minimal (we shall assume smoothness of their graph instead). Under some mild
assumptions, we shall characterize the existence of a desingularizing function
$\Psi$ that desingularizes the asymmetric modulus of the coderivative
(\emph{c.f.} Corollary~\ref{cor:des-as}). We give below our setting.

\subsection{Assumptions, setting}
\label{ss:31}
Let $S:\mathbb{R}\rightrightarrows\mathbb{R}^{n}$ be a multivalued map with
closed graph $\mathcal{S}$.

\begin{definition}
[smooth sweeping process]\label{smooth sweeping}We say that $S$ is a
\emph{smooth} sweeping process if either\smallskip\newline-- $\mathcal{S}$ is
a closed connected $\mathcal{C}^{1}$-smooth submanifold of $\mathbb{R}^{n+1}$
of dimension at most $n$ ; or\smallskip\newline-- $\mathcal{S}$ is a connected
smooth manifold of full dimension with boundary and $\partial\mathcal{S}$ is a
$\mathcal{C}^{1}$-smooth manifold of dimension $n$.
\end{definition}

\noindent It is clear that the above assumption is satisfied if $S$ is a
sweeping process associated to a $\mathcal{C}^{1}$-smooth function $f$
(\emph{c.f.} Example~\ref{example 1} or Remark~\ref{rem-1}). As a consequence
of this assumption we have the following result, which compares the modulus of $D^{\ast}S$
versus its asymmetric modulus. 

\begin{lemma}
\label{lem-compare-modula}  Let $S:\mathbb{R}\rightrightarrows \mathbb{R}^{n}$
be a smooth sweeping process and $(t,x)\in \mathcal{S}$. If either
\begin{equation*}
(a) \,\,\,\mathcal{S\;}\text{is a smooth manifold} \qquad \text{or} \qquad (b)\,\,\,
\|D^{\ast }S(t,x)|^{+}>0
\end{equation*}
then we have
\begin{equation*}
\Vert D^{\ast }S(t,x)|^{+}=\Vert D^{\ast }S(t,x)\Vert ^{+}.
\end{equation*}
\end{lemma}

\noindent\textbf{Proof. }If $\mathcal{S}$ is a smooth submanifold of
$\mathbb{R}^{n+1}$, the requested equality holds true for every $(t,x)\in
\mathcal{S}$ as a consequence of the fact that the limiting normal cone at any
point coincides with the normal space of the manifold at the same point. On
the other hand, if $\mathcal{S}$ is a manifold of full dimension with boundary
such that $\partial\mathcal{S}$ is also a smooth manifold, then the normal
cone $N_{\mathcal{S}}(t,x)$ is either $\{0\}$ or a ray generated by an outer
pointing normal vector $(s,y)$ of $\mathcal{S}$ at $(t,x)$. The conclusion
follows. \hfill$\Box$

\bigskip

\noindent Connectedness of $\mathcal{S}$ yields that $\mathrm{dom}
(S)$ is an interval (possibly unbounded). We shall use the following notation:
\[
\mathrm{T}=\sup(\mathrm{dom}(S)).\quad\text{(Notice that } T\in\mathbb{R}\cup\{+\infty\}.)
\]
We also define the multivalued map $H_{S}:\mathbb{R}\rightrightarrows
\mathbb{R}^{n+1}$ by
\[
H_{S}(t):=\partial\mathcal{S}\cap(\{t\}\times\mathbb{R}^{n}),\quad\text{for
all }t\in\mathbb{R}.
\]

\begin{ass*}
We say that $S$ satisfies the:

\begin{enumerate}
\item[$(A1)$] \emph{existence assumption} if for every $(t,x)\in\mathcal{S}$ with ${\Vert D^{\ast}S(t,x)|^{+}<+\infty}
$, there exist $\delta_x>0$ and at least one orbit $\gamma_{x}\in
AC(S;[t,t+\delta_{x}))\ $such that $\gamma_{x}
(t)=x$.

\item[$(A2)$] \emph{upper regular assumption} at $\overline{t}\in
\mathrm{dom}(S)$ with $\overline{t}<\mathrm{T}$, if there exists $\delta>0$
such that $\varphi^{\uparrow}(t)<+\infty$ for all $t\in(\overline{t}
,\overline{t}+\delta)$.

\item[$(A3)$] \emph{continuity assumption }at $\overline{t}\in\mathrm{dom}(S)$
with $\overline{t}<\mathrm{T}$, if there exists $\delta>0$ such that the
multivalued map $H_{S}$ is continuous for the Pompeiu-Hausdorff metric on
$(\overline{t},\overline{t}+\delta)$ (it may be discontinuous at
$\overline{t}$).
\end{enumerate}
\end{ass*}

\noindent Let us make some comments about the above assumptions: \smallskip\newline
Assumption (A1) ensures the existence of orbits issued from any non-critical
point. This assumption is satisfied if the sweeping process is defined via
(\ref{eq:Sf}) where $f$ is a $\mathcal{C}^{1,1}$-smooth function, since in
this case the existence of gradient orbits $\dot{\gamma}=-\nabla f(\gamma)$ is
guaranteed, and these orbits are also orbits for the sweeping process $S_{f}$ up to a suitable reparametrization, see Remark~\ref{remark 9}.
Assumption (A1) is also fulfilled if $S$ is a definable sweeping process, see
\cite[Section 6]{DD} or \cite{GRibarska}. In the general case, classical
existence results go back to the seminal work of J.J. Moreau \cite{M} for
convex-valued multifunctions which are Lipschitz continuous under the
Hausdorff-Pompieu metric. Since then, several extensions have been obtained,
see \cite{CG,CK,JV sweep} and references therein.\smallskip\newline Assumption
(A2) is automatically satisfied in the definable case, since in this case the
set of critical values is finite. In the general case, this assumption is
analogous to the hypothesis made in \cite[Section~3.3]{BDLM} that the critical
values of $f$ are upper isolated (see also statement of
Theorem~\ref{theorem BDLM}).\smallskip\newline Assumption (A3) is the more
restrictive, although it is natural in our setting. It is sa\-tisfied for
the sweeping process $S_{f}$ defined in (\ref{eq:Sf}) whenever $f$ is convex
or quasiconvex. In general, a\ smooth multivalued map $t\rightrightarrows
S(t)$ is not necessarily monotone in the sense of set-inclusion and the sets
$S(t)$ are not assumed convex (or of the same homology), therefore (A3) is
required to guarantee a control on the behavior of the boundaries. In
particular, the following result holds. (For the definitions of outer and inner semicontinuity 
of a multifunction the reader is referred to~\cite[\S 5]{RW}.)

\begin{proposition}
\label{Prop_seb}Let $S:\mathbb{R}\rightrightarrows\mathbb{R}^{n}$ be a smooth
sweeping process with bounded values and $a,b\in\mathbb{R}$ such that
$(a,b)\subset\mathrm{dom\,}(S)$. If $H_{S}$ is continuous on $(a,b)$, then $S$
is also continuous on $(a,b)$.
\end{proposition}

\noindent\textbf{Proof.} Let $I$ be a nontrivial interval contained in
a compact subset of $(a,b)$. It is sufficient to prove that $S$ is continuous on $I$. Since $\mathcal{S}\subset\mathbb{R}^{n+1}$ is closed and $S(t)=\mathcal{S}\cap (\{t\}\times \mathbb{R}^{n}),$
for every $t\in \mathbb{R}$, the map $S$ has closed (therefore, compact)
values and $S$ is outer semicontinuous. Let us assume, towards
a contradiction, that $S$ is not continuous on $I$, that is, there exists $\bar{t}\in I$ such that $S$ is not inner semicontinuous at 
$\bar{t}$. We deduce that there exist $\overline{x}\in S(\overline{t})$, $\varepsilon >0$ and a
sequence $\{t_{k}\}_{k}\subset \mathrm{dom\,}(S)$, converging to $\overline{t}$, such that
\begin{equation*}
d(\overline{x},S(t_{k}))\geq \varepsilon ,\quad \text{for all }k\in \mathbb{N}.
\end{equation*}
The above easily yields that $(\bar{t},\bar{x})\in \mathcal{S}\setminus 
\mathrm{int\,}(\mathcal{S})$, that is, $(\overline{t},\overline{x})\in
\partial \mathcal{S}$. However, since
\begin{equation*}
\mathcal{S}\,\cap\,\big(\,\{t_{k}\}\times B\left( \overline{x},\varepsilon \right) \big)
=\emptyset ,
\end{equation*}
this contradicts the continuity of $H_{S}$ at $\overline{t}$. \hfill $\Box $
\bigskip

\begin{remark}
In general, the converse of Proposition~\ref{Prop_seb} is not true. To see
this, set 
\begin{equation*}
\mathcal{S}:=\left( \mathbb{R}\times \lbrack -2,2]\right) \setminus \left\{
(t,x)\in \mathbb{R}^{2}:\,(t-1)^{2}+x^{2}\leq 1\right\} 
\end{equation*}%
and consider the sweeping process $S:\mathbb{R}\rightrightarrows \mathbb{R}$
defined by 
\begin{equation*}
S(t)=\mathcal{S\,}\cap \,\left( \{t\}\times \mathbb{R}^{2}\right) .
\end{equation*}%
It follows easily that $S$ is a smooth sweeping process. Moreover, $S$ is
continuous at every $t\in \mathbb{R}$, but $H_{S}$ is discontinuous at $0$.
\end{remark}

\subsection{Theorem A (characterizations via continuous dynamics)}
\label{ss:thmA}
Before we proceed, let us set
\[
\mathcal{T}:=\{t\in\mathrm{dom}(S):\,\text{(A2)--(A3) are fulfilled at }t\}.
\]
Observe that, if $t\in\mathcal{T}$, then there is $\delta>0$ such that
$[t,t+\delta)\subset\mathcal{T}$. \smallskip\newline We are now ready to state
the main result of this work. The proof will be given in Section~\ref{ss: proof of A}.

\begin{theoremA*}
Let $S:\mathbb{R}\rightrightarrows\mathbb{R}^{n}$ be a smooth sweeping process
with bounded values that satisfies~(A1). Let
$a\in\mathcal{T}$ (typically a critical value for $D^{\ast}S$). \smallskip\newline
The following assertions are equivalent:

\begin{enumerate}
\item[a)] \textbf{(Desingularization of the coderivative)} There exist $b>a$,
$\rho>0$ and a homeomorphism ${\Psi:[0,\rho]\rightarrow\lbrack a,b]}$, which
is a $\mathcal{C}^{1}$-diffeomorphism between $(0,\rho)$ and $(a,b)$ with
$\Psi^{\prime}(r)>0$ for every $r\in(0,\rho)$, such that:
\[
\Vert D^{\ast}(S\circ\Psi)(r,x)|^{+}\leq1,\quad~\text{for all }r\in
(0,\rho),~\text{for all }x\in S(\Psi(r)).
\]

\item[b)] \textbf{(Uniform length control for the absolutely continuous
orbits)} There exist $b>a$ and an increasing continuous function
$\sigma:[a,b]\mapsto\mathbb{R}^{+}$ with $\sigma(a)=0$ such that for every $a\leq t_{1}<t_{2}\leq
b$ and $\gamma\in\mathcal{AC}(S,[t_{1},t_{2}])$ we have:
\[
\ell(\gamma)\leq\sigma(t_{2})-\sigma(t_{1}).
\]

\item[c)] \textbf{(Length bound for the piecewise absolutely continuous
orbits)} There exist $b>a$ and $M>0$ such that for every $\gamma
\in\mathcal{PAC}(S,[a,b])$ we have:
\[
\ell(\gamma)\leq M.
\]

\item[d)] \textbf{(Integrability of the talweg)} There exists $b>a$ such that
\[
\int_{a}^{b}{\varphi^{\uparrow}}(t)<\infty.
\]

\end{enumerate}
\end{theoremA*}

\subsection{Theorem B (characterizations via discrete dynamics)}
\label{ss:thmB}
We first need the following definition.

\begin{definition}
[piecewise catching-up sequence]\label{discrete solution} Let $S:\mathbb{R}
\rightrightarrows\mathbb{R}^{n}$ be a multivalued map with closed
values.\smallskip\newline(i). A (finite or infinite) sequence $\{(t_{i}
,x_{i})\}_{i\geq0}\subset\mathcal{S}$ is called a \emph{catching-up sequence}
for $S\ $if $\{t_{i}\}_{i\geq0}$ is strictly increasing and
\[
x_{i+1}\in\mathrm{Proj}_{S(t_{i+1})}(x_{i}),\quad\text{for }i\geq0.
\]
(ii). A (finite or infinite) sequence of the form
\[
(t_{0}^{0},Y_{0}^{0}),(t_{1}^{0},Y_{1}^{0}),\ldots,(t_{k_{0}}^{0},Y_{k_{0}
}^{0}),(t_{0}^{1},Y_{0}^{1}),(t_{1}^{1},Y_{1}^{1}),\ldots,(t_{k_{1}}
^{1},Y_{k_{1}}^{1}),\ldots
\]
is called a \emph{piecewise catching-up sequence} for $S$ if for every
$j\geq0$
\[
\{(t_{i}^{j},Y_{i}^{j})\}_{i=0}^{k_{j}}\subset\mathcal{S}\text{ \ is a
catching-up sequence for }  S\,\,\,\text{and}\,\,t_{k_{j}}^{j}=t_{0}^{j+1}.
\]

\end{definition}

Now we are ready to state our second result which complements Theorem~A.

\begin{theoremB*}
\label{theorem B} The statements $(a)$-$(d)$ of Theorem~A are also equivalent to the following:

\begin{enumerate}
\item[e)] \textbf{(Uniform control of catching-up sequences)} There exist
$b>a$ and a continuous increasing function $\sigma:[a,b)\rightarrow
\lbrack0,\infty)$, with $\sigma(a)=0$, such that for every catching-up
sequence $\{(t_{i},x_{i})\}_{i\geq0}\subset\mathcal{S}$  with $\{t_{i}
\}_{i\geq 0}\in(a,b)$, and every $k\geq1$ we have
\begin{equation}\label{ortega}
\sum_{i=0}^{k}\Vert x_{i+1}-x_{i}\Vert\leq\sigma(t_{k})-\sigma(t_{0}).
\end{equation}

\item[f)] \textbf{(Length bound for piecewise catching-up sequences)} There
exist $b>a$ and $C>0$ such that for any piecewise catching-up sequence
\[
\left\{  (t_{i}^{j},Y_{i}^{j}):\;j\geq0,\;i\in\{0,\ldots,k_{j}\}\right\}
\]
with
\[
a<t_{0}^{0}<t_{1}^{0}<\ldots<t_{k_{0}}^{0}=t_{0}^{1}<t_{1}^{1}<\ldots<b
\]
we have:
\[
\sum_{j\geq0}\sum_{i=0}^{k_{j}}\Vert Y_{i+1}^{j}-Y_{i}^{j}\Vert\leq C.
\]

\end{enumerate}
\end{theoremB*}

\section{Proofs}

\label{Sec:Proofs}In this section we give proofs to our two main results,
Theorem~A (Subsection~\ref{ss: proof of A}) and Theorem~B
(Subsection~\ref{ss: proof of B}). To do so, we shall need some auxiliary
results (Subsection~\ref{ss:aux}) and a new notion of oriented calmness
(Subsection~\ref{ss:calm}).

\subsection{Auxiliary results}

\label{ss:aux}The first result concerns continuity of the moduli maps. It is
based on the fact that the normal space mapping of the smooth manifold is
continuous (in the Grasmannian metric). The details are left to the reader.

\begin{lemma}
[continuity of the (oriented) modulus on $\partial\mathcal{S}$]
\label{continuity of the coderivative} Let $S:\mathbb{R}\rightrightarrows
\mathbb{R}^{n}$ be a smooth sweeping process. Then, the functions
\[
{(t,x)\mapsto\Vert D^{\ast}S(t,x)|^{+}\qquad}\text{and}{\qquad(t,x)\mapsto
\Vert D^{\ast}S(t,x)\Vert^{+}}
\]
are continuous on $\partial\mathcal{S}$ for the usual topology on
$\mathbb{R}\cup\{+\infty\}$.
\end{lemma}

\noindent The second result asserts continuity of the (oriented) talweg function.
Let us recall from Subsection~\ref{ss:31} that the multivalued function 
$H_S:\mathbb{R}\rightrightarrows\mathbb{R}^n$ is defined by $H_{S}(t):=\partial\mathcal{S}\cap(\{t\}\times\mathbb{R}^{n})$, for
all $t\in\mathbb{R}$.

\begin{lemma}
[continuity of the (oriented) talweg function]\label{continuity of the talweg}
Let $S:\mathbb{R}\rightrightarrows\mathbb{R}^{n}$ be a smooth sweeping process
such that $S(t)$ is bounded for all $t\in\mathbb{R}$. Let $[a,b]\subset
\mathrm{dom}(S)$ such that $H_S$ is continuous for the
Pompeiu-Hausdorff metric on $[a,b]$. Then the talweg functions $\varphi
^{\uparrow}$ and $\varphi$ are continuous on $[a,b]$, where the image space
$\mathbb{R}\cup\{+\infty\}$ is considered with its usual topology.
\end{lemma}

\noindent\textbf{Proof. }Set $K:= H_S([a,b])$, which is a compact set.  
Since
\[
\varphi^{\uparrow}(t)=\,\underset{x\in H_S(t)}{\max}{\Vert D^{\ast}S(t,x)|^{+}
}\qquad\text{(respectively, \ }\varphi(t)=\,\underset{x\in H_S(t)}{\max}{\Vert
D^{\ast}S(t,x)||^{+}}\text{).}
\]
the result follows from Lemma~\ref{continuity of the coderivative}.\hfill
$\Box$

\bigskip

\begin{proposition}
[diffeomorphic rescaling of time]\label{modifying time} Let $S:\mathbb{R}
\rightrightarrows\mathbb{R}^{n}$ be a multivalued map and $\gamma
\in\mathcal{AC}(S,(a,b))$. If $\Psi:(0,\rho)\rightarrow(a,b)$ is a $\mathcal{C}
^{1}$-smooth diffeomorphism such that $\Psi^{\prime}(r)>0$ for all $r\in
(0,\rho)$, then $\tilde{\gamma}=\gamma\circ\Psi$ is an orbit of the sweeping
process defined by $\tilde{S}:=S\circ\Psi$, that is, $\tilde{\gamma}
\in\mathcal{AC}(\tilde{S},(0,\rho)).$
\end{proposition}

\noindent\textbf{Proof. }It is straighforward that $\tilde{\gamma}=\gamma
\circ\Psi$ is an absolutely continuous curve. Since $\Psi$ is a bi-Lipschitz
homeomorphism on each compact interval contained in $(0,\rho)$ we deduce that for
any null subset $A$ of $(a,b)$ the set $\Psi^{-1}(A)\ $is also null (with
respect to the Lebesgue measure). If $\mathcal{I}$ be the points of
differentiability of $\gamma$ for which~\eqref{sweepingequation} holds, it
follows that $\mathcal{J}:=\Psi^{-1}((a,b)\setminus\mathcal{I})$ is a null set
and for every $r\in(0,\rho)\setminus\mathcal{J}$ it holds:
\[
\tilde{\gamma}^{\prime}(r)=(\gamma\circ\Psi)^{\prime}(r)=\gamma^{\prime}
(\Psi(r))\Psi^{\prime}(r)\in N_{S(\Psi(r))}(\gamma(\Psi(r))),
\]
yielding that $\tilde{\gamma}$ is an orbit solution of the sweeping process
defined by $S\circ\Psi$.\hfill$\Box$

\bigskip

\noindent In the sequel, given a curve $\gamma:I\rightarrow\mathbb{R}^{n}$ we define its
lifting $\zeta:I\rightarrow\mathbb{R}^{n+1}$ by
\[
\zeta(t)=(t,\gamma(t)),{\qquad t\in I.}
\]

\begin{proposition}
[geometric facts]\label{geometric facts} Let $S:\mathbb{R}\rightrightarrows
\mathbb{R}^{n}$ be a smooth sweeping process. \\  Fix $\overline{t}\in
\mathrm{dom}(S)\setminus\{\mathrm{T}\}$ and $\bar{x}\in S(\overline
{t})$. Then:

\begin{enumerate}
\item[$a)$] If there is $\delta>0$ such that $\bar{x}\in S(t)$, for all $t\in(\overline{t},\overline
{t}+\delta)$, then $\alpha\leq0$ for all $(\alpha,u)\in
N_{\mathcal{S}}(\overline{t},\bar{x})$.

\item[$b)$] If $\Vert D^{\ast}S(\overline{t},\bar{x})|^{+}>0$, then for any
$\tau>\overline{t}$ and $\gamma\in\mathcal{AC}(S,[\overline{t},\tau))$ with
$\gamma(\overline{t})=\bar{x}$, there exists $\delta>0$ such that
\[
\zeta(t):=(t,\gamma(t))\in\partial\mathcal{S},{\qquad}\text{for all }
t\in\lbrack\overline{t},\overline{t}+\delta).
\]

\item[$c)$] If $\mathrm{int}(\mathcal{S})$ is nonempty and $N_{\mathcal{S}
}(\overline{t},\bar{x})=\mathbb{R}_{+}(\alpha,u)$ with $\alpha<0$, then there
is $\delta>0$ such that $\bar{x}\in S(t)$ for all $t\in\lbrack\overline
{t},\overline{t}+\delta)$.
\end{enumerate}
\end{proposition}

\noindent\textbf{Proof. }$(a).$\textbf{ }If $(\overline{t},\bar{x}
)\in\mathrm{int}$\thinspace$(\mathcal{S})$ then $N_{\mathcal{S}}(\bar{t}
,\bar{x})=\{(0,0)\}$ and the conclusion follows trivially. In the case when $(\overline{t},\overline{x})\in \partial\mathcal{S}$, since $\partial\mathcal{S}$ is a smooth manifold, the limiting normal
cone $N_{\mathcal{S}}(\overline{t},\bar{x})$ is equal to the Fr\'{e}chet
normal cone and is contained in the normal space of $\partial\mathcal{S}$ at
$(\overline{t},\bar{x})$. Therefore, for any $(\alpha,u)\in N_{\mathcal{S}
}(\bar{t},\bar{x})$ and $t\in(\overline{t},\overline{t}+\delta)$, we have
$(t,\bar{x})\in\mathcal{S}$ and
\[
\limsup_{t\searrow\overline{t}}\dfrac{\langle(\alpha,u),(t-\overline{t},\bar
{x}-\bar{x})\rangle}{\Vert(t-\overline{t},\bar{x}-\bar{x})\Vert}=\alpha\leq0.
\]

$(b).$ Let $\tau >\overline{t}$ and $\gamma \in \mathcal{AC}(S,[\overline{t},\tau ))$ with $\gamma (\overline{t})=\bar{x}$ and assume $\Vert D^{\ast }S(\overline{t},\bar{x})|^{+}>0.$ Since $(t,y)\mapsto \Vert D^{\ast }S(t,y)|^{+}
$ is continuous on $\partial \mathcal{S}$ (Lemma~\ref{continuity of the
coderivative}), there exists a neighborhood $\mathcal{V}$ of $(\overline{t},
\bar{x})$ such that for all $(t,y)\in \mathcal{V}\cap \partial \mathcal{S}$
we have $\Vert D^{\ast }S(t,y)|^{+}>0$. Therefore, there is $\delta >0$ such that $\Vert D^{\ast }S(\zeta (t))|^{+}>0$
and consequently, $\zeta (t)\in \partial \mathcal{S}$ for all $t\in \lbrack 
\overline{t},\overline{t}+\delta )$.\smallskip 

$(c).$ It follows from our assumption that $\mathrm{dim\,}(\partial \mathcal{S})=n$ and $(\alpha ,u)$ is a nonzero outer normal vector of $\mathcal{S}$
at $(\overline{t},\bar{x})$. 
Without loss of generality, let us assume that $(\alpha,u)$ is a unit vector.
Since  $\mathrm{int}(\mathcal{S})\neq \emptyset
,$ we deduce that ${(\overline{t},\bar{x})-}\lambda (\alpha ,u)\in \mathcal{S}$ for all $\lambda >0$ sufficiently small. Let us assume, reasoning towards a contradiction, that there exists a decreasing sequence ${\{t_{k}\}_{k}\subset \mathbb{R}}$ converging to $\overline{t}$ such that $\bar{x}\notin
S(t_{k})$, for all $k\in \mathbb{N}$.  Let us now take a decreasing sequence 
$\{\lambda _{k}\}_{k}\subseteq \mathbb{R}^{+}$ that converges to $0$ and
satisfies $(\overline{t},\bar{x})-\lambda _{k}(\alpha ,u)\in \mathcal{S}$
for all $k$. Let $\boldsymbol{z}_{k}\in\mathbb{R}^{n+1}$ be any vector such that 
\begin{equation*}
\boldsymbol{z}_{k}\in \partial \mathcal{S}\,\bigcap \,\, \left [\, (t_{k},\bar{x}),\,(\overline{t}-\lambda _{k}\alpha ,\bar{x}-\lambda _{k}u)\,\right ]\,,
\end{equation*}
where $\left [\,(t_{k},\bar{x}),\,(\overline{t}-\lambda _{k}\alpha ,\bar{x}-\lambda_{k}u)\,\right ]$ stands for the line segment joining the points $(t_{k},\bar{x})$ and $(\overline{t}-\lambda _{k}\alpha ,\bar{x}-\lambda _{k}u)$. 
It follows easily that $\{\mathbf{z}_{k}\}_{k}$ converges to $(\overline{t},\bar{x})$ and that 
\[\left\langle \dfrac{\mathbf{z}_k -(\overline{t},\overline{x})}{\| \mathbf{z}_k -(\overline{t},\overline{x})\|}, (\alpha,u)\right\rangle \,\leq\, \left\langle \,(1,0),(\alpha,u)\,\right\rangle = \alpha.\]
Let $d$ be any accumulation point of the sequence $(\boldsymbol{z}_{k}-(\overline{t},x))/\Vert \boldsymbol{z}_{k}-(\overline{t},x)\Vert $. Then, $\boldsymbol{d}$ belongs to the Bouligand tangent cone of $\partial \mathcal{S}$, which
coincides with the tangent space of $\mathcal{S}$ at the same point. Therefore $\boldsymbol{d}$ should be orthogonal to the normal vector $(\alpha,u)$. However, $\langle \boldsymbol{d}, (\alpha,u)\rangle\leq \alpha <0$, which leads to a
contradiction.  \hfill $\Box$

\bigskip

\noindent The following lemma is crucial in the proof of our main theorem
since it relates the value of the coderivate with the velocity of the orbit
of the sweeping process. The proof follows closely the proof of \cite[Theorem
4.1]{DD} where a similar result has been established for the usual modulus $%
\Vert D^{\ast }S(t,\gamma (t))\Vert ^{+}$.

\begin{lemma}
\label{boundedvelocity} Let $S:\mathbb{R}\rightrightarrows\mathbb{R}^{n}$ be a
smooth sweeping process and $\gamma\in\mathcal{AC}(S,[a,b))$. Then,
\[
\Vert\dot{\gamma}(t)\Vert=\Vert D^{\ast}S(t,\gamma(t))|^{+},
\]
for all $t\in\lbrack a,b)$ such that $-\dot{\gamma}(t)\in N_{S(t)}(\gamma(t))$
and $\Vert D^{\ast}S(t,\gamma(t))|^{+}$ is finite.
\end{lemma}

\noindent\textbf{Proof. }Let $t\in\lbrack a,b)$ be a point of
differentiability of $\gamma$ such that $-\dot{\gamma}(t)\in N_{S(t)}
(\gamma(t))$ and that $\Vert D^{\ast}S(t,\gamma(t))|^{+}$ is finite.\bigskip

\noindent \emph{First case}: $\dot{\gamma}(t)=0$.\smallskip \newline
If $\zeta(t):=(t,\gamma(t))\in
\mathrm{int}$\thinspace$(\mathcal{S})$, the desired equality holds trivially,
while if $\zeta(t)\in\partial \mathcal{S}$, then ${\dot{\zeta}(t)=(1,0)}$
belongs to the tangent space of $\partial\mathcal{S}$ at $\zeta(t)$. Since $S$
is a smooth sweeping process, the normal cone $N_{\mathcal{S}}(\zeta(t))$ is
contained in the normal space of $\partial\mathcal{S}$ at $\zeta(t)$.
Therefore,
\[
\langle(1,0),N_{\mathcal{S}}(\zeta(t))\rangle=\{0\}.
\]
Hence, if $(\alpha,u)\in N_{\mathcal{S}}(\zeta(t))$, then $\alpha=0.$ Thus,
$\Vert D^{\ast}S(\zeta(t))|^{+}=0$.\bigskip

\noindent \emph{Second case}: $\dot{\gamma}(t)\neq0$. \smallskip\newline Then
$\zeta(t)\in\partial\mathcal{S}$ and  $\dot{\zeta}(t)$ belongs to the tangent
space of $\partial\mathcal{S}$ at $\zeta(t)$. As in the first case, we obtain
that
\[
\langle(1,\dot{\gamma}(t)),N_{\mathcal{S}}(\zeta(t))\rangle=\{0\}.
\]
Hence, for every $(\alpha,u)\in N_{\mathcal{S}}(\zeta(t))$ with $\Vert
u\Vert=1$ we have $\alpha+\langle\dot{\gamma}(t),u\rangle=0$. Thanks to
Cauchy-Schwartz inequality, we obtain
\[
\Vert\dot{\gamma}(t)\Vert\geq\Vert D^{\ast}S(\zeta(t)|^{+}.
\]
By Proposition~\ref{geometric facts} $(c)$, we can assume that
\[
\underset{\Vert u\Vert\leq1}{\sup}\big\{  a:\,a\in D^{\ast}S(t,x)(u)\big\}
\geq0.
\]
Setting $H=\{t\}\times\mathbb{R}^{n}$ we have $\{t\}\times S(t)=H\cap
\mathcal{S}.$ Due to the fact that $-\dot{\gamma}(t)\in N_{S(t)}(\gamma(t))$, we have:
\[
(1,-\dot{\gamma}(t))\in N_{\{t\}\times S(t)}(\zeta(t)).
\]
In addition, since $S$ is a smooth sweeping process and ${\Vert D^{\ast}S(t,\gamma(t))|^{+}<\infty}$, we have that
$(t,0)\in N_{\mathcal{S}}(t,\gamma(t))$ only if $t=0$. Hence, applying the
calculus rule \cite[Theorem~6.42]{RW}, we get
\[
N_{H\cap\mathcal{S}}(\zeta(t))\subset N_{H}(\zeta(t))+N_{\mathcal{S}}
(\zeta(t))=\mathbb{R}\times\{0\}+N_{\mathcal{S}}(\zeta(t)).
\]
Therefore, the inclusion $(\lambda,-\dot{\gamma}(t))\in N_{\mathcal{S}}
(\zeta(t))$ holds for some $\lambda\in\mathbb{R}$. By orthogonality between normal and tangent vectors, we
get that:
\[
\langle(\lambda,-\dot{\gamma}(t)),(1,\dot{\gamma}(t))\rangle=0.
\]
and thus $\lambda=\Vert\dot{\gamma}(t)\Vert^{2}$. After
normalization, we obtain:
\[
\left(  \Vert\dot{\gamma}(t)\Vert,-\dfrac{\dot{\gamma}(t)}{\Vert\dot{\gamma
}(t)\Vert}\right)  \in N_{\mathcal{S}}(\zeta(t)),
\]
which readily yields $\Vert D^{\ast}S(t,\gamma(t))|^{+}\geq\Vert\dot{\gamma
}(t)\Vert$, as claimed.\hfill$\Box$

\bigskip

Let us finally quote the following result, which is a restatement of (and can
be proved in the same way as) \cite[Proposition~27]{BDLM}.

\begin{proposition}[concatenation]
\label{construction} Let $b>a$ and $\Gamma$ be a collection of absolutely
continuous curves $\gamma$ defined in some nontrivial interval $J\subset(a,b)$
with values in $\mathbb{R}^{n}$. Assume that for each $t\in(a,b)$ there exist
$\varepsilon_{t}>0$ and $\gamma_{t}\in\Gamma$ with $\mathrm{dom\,}
(\gamma_t)=[t,t+\varepsilon_{t})$. Then there exist a countable partition
$\{I_{n}\}_{n\in\mathbb{N}}$ of $(a,b)$ into intervals $I_{n}$ of nonempty
interior and a piecewise absolutely continuous curve $\gamma: (a,b) \rightarrow
\mathbb{R}$ such that for each $n\in\mathbb{N}$, there is $\gamma^n\in \Gamma$ such that $\gamma=\gamma^n$ on $I_n$.
\end{proposition}

We are now ready to prove our main result.

\subsection{Proof of Theorem A}
\label{ss: proof of A}
 We prove $(a)\Rightarrow (b)\Rightarrow (c)\Rightarrow
(d)\Rightarrow (a).$\medskip \newline
$\bm{a) \to b)}:$ Let $\Psi:[0,\rho]\rightarrow\lbrack a,b]$ be
given by $(a)$. 
Let $\gamma\in \mathcal{AC}(S, [t_1,t_2))$ with $[t_1,t_2)\subset [a,b)$. Since $\Psi$ is a $\mathcal{C}^{1}
$-smooth function, $\partial\mathrm{gph}((S\circ\Psi)|_{(0,\rho)})$ is a smooth
manifold. By Proposition~\ref{modifying time}, $\gamma\circ\Psi \in \mathcal{AC}(S\circ\Psi, \Psi^{-1}([t_1,t_2)))$.
Applying Lemma~\ref{boundedvelocity}, we deduce that
\[
\left\vert \dfrac{d\left(  \gamma\circ\Psi\right)  }{dr}(r)\right\vert =\Vert
D^{\ast}(S\circ\Psi)(r,\gamma(\Psi(r)))|^{+}\leq1,\quad\forall_{a.e}
~r\in(a,d).
\]
Since $\Psi$ is increasing and smooth, by change of variables we obtain:
\begin{align*}
\int_{t_{1}}^{t_{2}}\Vert\dot{\gamma}(\tau)\Vert d\tau=\int_{\Psi^{-1}(t_{1}
)}^{\Psi^{-1}(t_{2})}\left\Vert \dot{\gamma}(\Psi(r))\right\Vert \dot{\Psi
}(r)dr &  =\int_{\Psi^{-1}(t_{1})}^{\Psi^{-1}(t_{2})}\left\Vert \dfrac
{d\left(  \gamma\circ\Psi\right)  }{dr}(r)\right\Vert dr\\
&  \leq\int_{\Psi^{-1}(t_{1})}^{\Psi^{-1}(t_{2})}dr=\Psi^{-1}(t_{2})-\Psi
^{-1}(t_{1}).
\end{align*}
Therefore $(b)$ is satisfied by setting $\sigma:=\Psi^{-1}.$ \medskip\newline
$\bm{b) \to c)}:$ Since $\sigma$ is an increasing function and $\sigma
(a)=0,$ statement $(c)$ follows by setting $M:=\sigma(b).$ \medskip\newline
$\bm{c) \to d)}:$ Let $b>a$ and let $M>0$ given by statement $c)$. Let
${\varphi^{\uparrow}:(a,b)\rightarrow\mathbb{R}\cup\{+\infty\}}$ be the
oriented talweg function of $S$ and let us assume, towards a contradiction,
that for any $c\in(a,b)$ the function ${\varphi^{\uparrow}}$ is not integrable
on $(a,c)$. By Lemma~\ref{continuity of the coderivative}, the function
$(t,x)\mapsto\Vert D^{\ast}S(t,x)|^{+}$ is continuous on $\partial\mathcal{S}
$. By assumptions (A2)--(A3), shrinking $b$ if necessary, we may assume that
${\varphi^{\uparrow}}(t)<\infty$ for all $t\in(a,b)$ and that the multivalued
map $t\rightrightarrows H_{S}(t)$ is continuous on $(a,b)$. 
By
Lemma~\ref{continuity of the talweg}, $\varphi^{\uparrow}$ is continuous on
$(a,b)$.\smallskip\newline By Lemma~\ref{boundedvelocity}, if $J$ is a
nontrivial interval of $(a,b)$ then for any $\gamma\in\mathcal{AC}(S,J)$ we
have $\Vert\dot{\gamma}(t)\Vert=\Vert D^{\ast}S(t,\gamma(t))|^{+}$ for almost
every $t\in J$. Let $k\in\mathbb{N}$ and $t\in(a,b)$ and define a curve
$\gamma_{t}^{k}$ as follows:

\begin{itemize}
\item If ${\varphi^{\uparrow}}(t)=0$, take $\gamma_{t}^{k}\in\mathcal{AC}
(S,[t,\tau))$ be any curve such that ${\tau-t<1/k}$.

\item If ${\varphi^{\uparrow}}(t)>0$, since $H_{S}(t)$ is compact, there
exists $x\in S(t)$ such that $\Vert D^{\ast}S(t,x)|^{+}=\varphi^{\uparrow}
(t)$. Thanks to assumption (A1) and Lemma~\ref{continuity of the talweg}, we
can take $\gamma_{t}^{k}\in\mathcal{AC}(S,[t,\tau)),$ for some $\tau>t,$ such
that $\gamma(t)=x$ and
\[
\Vert\dot{\gamma}_{t}^{k}(s)\Vert>\frac{k-1}{k}\varphi^{\uparrow}
(s),\quad\text{for almost every }s\in(t,\tau).
\]
\end{itemize}
Gluing together, thanks to Proposition~\ref{construction} (concatenation), we obtain
$\gamma^{k}\in\mathcal{PAC}(S,(a,b))$ such that for almost every $t\in(a,b)$
\[
\varphi^{\uparrow}(t)\geq\Vert\dot{\gamma}^{k}(t)\Vert\geq f_{k}(t):=
\begin{cases}
\phantom{xx}0,\medskip & \text{if }t\in A_{k}\\
\dfrac{k-1}{k}\varphi^{\uparrow}(t), & \text{if }t\in(a,b)\setminus A_{k}.
\end{cases}
\]
where $A=\{t\in(a,b):~\varphi^{\uparrow}(t)=0\}$ and $A_{k}=(a,b)\cap
(A+[0,1/k])$ for all $k\in\mathbb{N}$.\smallskip\newline The continuity of
$\varphi^{\uparrow}$ yields that $A$ is a closed set relatively to $(a,b)$.
Therefore, ${A=\cap_{k\in\mathbb{N}}A_{k}}$. Then, for all $t\in(a,b)$,
$f_{k}(t)\nearrow\varphi^{\uparrow}(t)$ as $k$ tends to infinity. Hence, by
the Monotone Convergence Theorem, $(\int_{a}^{b}f_{k})_{k}$ converges to
$\int_{a}^{b}\varphi^{\uparrow}$, which is infinity. Thus, there is
$K\in\mathbb{N}$ such that
\[
\int_{a}^{b}\Vert\dot{\gamma}^{K}(t)\Vert dt\geq\int_{a}^{b}f_{K}(t)dt>M,
\]
which contradicts statement $(c)$ since $\gamma^{K}\in\mathcal{PAC}
(S,(a,b))$.\medskip\newline
$\bm{d)\to a)}:$ Let us assume that the oriented talweg function
$\varphi^{\uparrow}$ is integrable on $[a,b]$ for some $b>a.$ As a consequence
of assumptions (A2) and (A3), shrinking $b$ if necessary, we may assume that
${\varphi^{\uparrow}}$ is continuous on $[a,b]$ and $\varphi^{\uparrow
}(t)<\infty$ for all $t\in(a,b]$. Let $\overline{\varphi}:=\max\{\varphi
^{\uparrow},1\}$ which is an integrable continuous majorant of ${\varphi^{\uparrow}
}$ and set
\[
\theta(t):=\int_{a}^{t}\overline{\varphi}(s)ds,\quad\text{for }t\in\lbrack
a,b].
\]
Since $\overline{\varphi}$ is positive and integrable on $[a,b]$, we set
$\rho:=\theta(b)$ and define ${\Psi:[0,\rho]\rightarrow\lbrack a,b]}$ as the
inverse function of $\theta$, that is, $\Psi(r)=\theta^{-1}(r)$. Since
$\theta^{\prime}(t)=\overline{\varphi}(t)\in\lbrack1,+\infty),$ for every
$t\in(a,b],$ it follows that $\Psi$ is $\mathcal{C}^{1}$-smooth on $(0,\rho)$,
with derivative
\[
{\Psi}^{\prime}(r)=\frac{1}{\overline{\varphi}(\Psi(r))}\leq1,\quad\text{for
all }r\in(0,\rho).
\]
Thus, $\Psi$ is a Lipschitz homeomorphism between $[0,\rho]$ and $[a,b]$.
Finally, using the chain rule for coderivatives \cite[Theorem 10.37]{RW}, we
deduce that
\[
\Vert D^{\ast}\left(  S\circ\Psi\right)  (r,x)|^{+}\leq\dfrac{\Vert D^{\ast
}S(\Psi(r),x)|^{+}}{\overline{\varphi}(\Psi(r))}\leq1,\quad\text{for all }
r\in(0,\rho).
\]
The proof is complete. \hfill$\Box$

\subsection{Oriented calmness}

\label{ss:calm}Before proceeding to the proof of Theorem~B, we need to
introduce the modulus of oriented calmness and establish a result analogous
to the Mordukhovich criterium for the oriented modulus of the coderivative.
Let us first recall that the \emph{Lipschitzian graphical modulus} of
$S:\mathbb{R}\rightrightarrows\mathbb{R}^{n}$ at $t$ for $x$ is defined by
\begin{align*}
\mathrm{Lip\,}S(t,x):=\inf\{\kappa>0|~ &  \exists\epsilon>0,~\delta
>0,~\text{such that }\\
&  S(t_{2})\cap B(x,\delta)\subset S(t_{1})+\kappa|t_{2}-t_{1}|B, \quad\text{for all
 } t_{1}, t_{2}\in(t-\epsilon,t+\epsilon)\},
\end{align*}
where $B$ stands for the open unit ball. \smallskip\newline
We recall that the multivalued function $S$ has the Aubin property at $t$ for $x$ if
and only if $\mathrm{Lip\,}S(t,x)<\infty$. More precisely, we have the
following (see \cite[Theorem 9.40]{RW}).

\begin{theorem}
\label{Mordukhovich}For every $(t,x)\in\mathcal{S}$ such that $\Vert D^{\ast
}S(t,x)\Vert^{+}<\infty$ it holds:
\[
\mathrm{Lip\,}S(t,x)=\Vert D^{\ast}S(t,x)\Vert^{+}.
\]

\end{theorem}

\noindent Motivated by the above, we introduce the following graphical modulus.

\begin{definition}
[oriented calm modulus]Let $S:\mathbb{R}\rightrightarrows\mathbb{R}^{n}$ be a
multivalued map and $(t,x)\in\mathcal{S}$. The oriented calm graphical
modulus, denoted by $\mathrm{calm}^{\uparrow}S$, at $t$ for $x$ is defined  by
\begin{align*}
\mathrm{calm}^{\uparrow}S(t,x):= &  \inf\{\kappa>0|~\exists\epsilon
>0,~\delta>0,~\text{such that }\\
&  S(t)\cap B(x,\delta)\subset S(t_{1})+\kappa|t_{1}-t|B~\text{for all }
t_{1}\in(t,t+\epsilon)\}.
\end{align*}

\end{definition}

\noindent Observe that, if $S$ is a single-valued function and $\mathrm{calm}^{\uparrow}S(t,x)<\infty$, then $S$ is calm at~$t$ to the right. More information on the notion of calmness for multivalued maps can be found in~\cite{HO} and references therein. 
We are now ready to give the oriented version of Theorem~\ref{Mordukhovich}.

\begin{proposition}
[oriented calm vs oriented modulus]\label{oriented criterion}Let
$S:\mathbb{R}\rightrightarrows\mathbb{R}^{n}$ be a smooth sweeping process,
$t\in\mathrm{dom}(S)\setminus\{\mathrm{T}\}$ and $x\in S(t)$ such that $\Vert
D^{\ast}S(t,x)|^{+}<+\infty$. Then
\[
\mathrm{calm}^{\uparrow}S(t,x)=\Vert D^{\ast}S(t,x)|^{+}.
\]

\end{proposition}

\noindent\textbf{Proof.} Let us first notice that $\mathrm{calm}^{\uparrow }S(t,x)\leq \mathrm{Lip\,}S(t,x)$. We consider two cases:
\medskip \newline
\textit{Case 1}: $\Vert D^{\ast }S(t,x)|^{+}=0.$\smallskip \newline
If $\Vert D^{\ast }S(t,x)\Vert ^{+}=0,$ then $\mathrm{calm}^{\uparrow
}S(t,x)=0.$ If $\Vert D^{\ast }S(t,x)\Vert ^{+}>0$, then, by Lemma~\ref{lem-compare-modula} $\mathcal{S}$ is a manifold of full dimension with boundary $\partial \mathcal{S}$ which is a smooth manifold of dimension $n$.
Let us assume by contradiction that $\mathrm{calm}^{\uparrow }S(t,x)>0$.
Then, for every $k\in \mathbb{N}$ such that $k^{-1}<\mathrm{calm}^{\uparrow
}S(t,x)$, there exists $y_{k}\in S(t)\cap B(x,1/k)$ such that 
\begin{equation*}
y_{k}\notin S(t_{k})+\left( \dfrac{t_{k}^{\prime }-t}{k}\right) B,~\text{for
some }t_{k}^{\prime }\in (t,~t+\dfrac{1}{k}).
\end{equation*}
Set $t_{k}:=\inf \{r\in (t,t+\frac{1}{k}):~y_{k}\notin S(r)\}$. It is clear
that $(t_{k},y_{k})\in \partial \mathcal{S}$ and that $y_{k}$ is not
right-locally stationary for $S$ at $t_{k}$. Thus, by Proposition~\ref{geometric facts} $(c)$, for every $k\in \mathbb{N}$ and $(\beta
_{k},v_{k})\in N_{\mathcal{S}}(t_{k},y_{k})$, we have $\beta _{k}\geq 0$.
Since $N_{\mathcal{S}}(t,x)$ is a ray and $\{(t_{k},y_{k})\}_{k}\rightarrow
(t,x)$, the continuity of unit outer normal vectors of $\mathcal{S}$ on $\partial \mathcal{S}$ ensures that $\beta \geq 0$ whenever ${(\beta ,v)\in N_{\mathcal{S}}(t,x)}$. This leads to the equality $\Vert D^{\ast }S(t,x)|^{+}=\Vert
D^{\ast }S(t,x)\Vert ^{+}$, which is a contradiction. Therefore, $\mathrm{calm}^{\uparrow }S(t,x)=0$.  \smallskip\newline 
\textit{Case 2}: $\Vert D^{\ast }S(t,x)|^{+}=\alpha>0.$ \medskip \newline
In this case, we deduce from Lemma~\ref{lem-compare-modula}$(b)$ that 
\begin{equation*}
\Vert D^{\ast }S(t,x)|^{+}=\Vert D^{\ast }S(t,x)\Vert ^{+}=\mathrm{Lip\,}S(t,x)\geq \mathrm{calm}^{\uparrow }S(t,x).
\end{equation*}
By Lemma~\ref{continuity of the coderivative} and compactness of the unit
ball of $\mathbb{R}^{n},$ there exists $u\in \mathbb{R}^{n}$ with $\Vert
u\Vert =1$ such that $(\alpha ,u)\in N_{\mathcal{S}}(t,x)$. Let $\{t_{k}\}_{k\geq 1}\subset \mathbb{R}$ be a decreasing sequence that
converges to $t$. Let $\{y_{k}\}_{k\geq 1}\subset \mathbb{R}^{n}$ be a
sequence that satisfies $y_{k}\in \mathrm{Proj}(x,S(t_{k}))$ for each $k\in 
\mathbb{N}$. By compactness of the unit sphere of $\mathbb{R}^{n+1}$, up to
a subsequence we deduce that 
\begin{equation*}
\lim_{k\rightarrow \infty }\dfrac{(t_{k}-t,y_{k}-x)}{\Vert
(t_{k}-t,y_{k}-x)\Vert }=(\beta ,v),
\end{equation*}
where $(\beta ,v)$ belongs to the tangent space of $\mathcal{S}$ at $(t,x)$
and $\beta \geq 0.$ Since $\mathcal{S}$ is a smooth sweeping process, it follows that 
\begin{equation*}
(\alpha ,u)\bot (\beta ,v)\qquad \text{yielding\qquad }\langle u,v\rangle
=-\alpha \beta .
\end{equation*}
Since $\mathrm{calm}^{\uparrow }S(t,x)\leq \Vert D^{\ast
}S(t,x)|^{+}<+\infty $, $\beta $ must be a strictly positive number.
Therefore 
\begin{equation*}
\lim_{k\rightarrow \infty }\dfrac{\Vert y_{k}-x\Vert }{t_{k}-t}=\dfrac{\Vert
v\Vert }{\beta }\geq \dfrac{|\langle u,v\rangle |}{\beta }=\alpha ,
\end{equation*}
implying that 
\begin{equation*}
\mathrm{calm}^{\uparrow }S(t,x)\geq \alpha =\Vert D^{\ast }S(t,x)|^{+}.
\end{equation*}

\noindent The proof is complete.\hfill $\Box $

\bigskip

\begin{lemma}[controlling excess of $S(t_{0})$]
\label{distance-projection} Let $S:\mathbb{R}\rightrightarrows \mathbb{R}^{n}
$ be a smooth sweeping process and $[t_{0},t_{1}]\subset \mathrm{dom}(S)$.
Then 
\begin{equation*}
\mathrm{ex\,}(S(t_{0}),S(t_{1})):=\,\underset{x\in S(t_{0})}{\sup }
d(x,S(t_{1}))\leq \left( \underset{t\in \lbrack t_{0},t_{1}]}{\sup }
\varphi ^{\uparrow }(t)\right) (t_{1}-t_{0})
\end{equation*}
and 
\begin{equation*}
\mathrm{dist}(S(t_{0}),S(t_{1}))\leq \left( \underset{t\in \lbrack
t_{0},t_{1}]}{\sup }\varphi (t)\right) (t_{1}-t_{0}).
\end{equation*}
\end{lemma}

\noindent \textbf{Proof.} Let us first notice that 
\begin{equation*}
K:=\,\underset{t\in \lbrack t_{0},t_{1}]}{\sup }
\varphi ^{\uparrow }(t)\,\geq \|D^{\ast}S(t,x)|^+=\mathrm{calm}^{\uparrow }S(t,x),\quad\text{for all } t\in [t_0,t_1] \,\,\text{and } x\in S(t)\,.
\end{equation*}
 If $K=\infty $, there is nothing to prove. Let $K<+\infty$ and assume, towards a contradiction, that for some $\delta>0$ we have
\begin{equation*}
\mathrm{ex\,}(S(t_{0}),S(t_{1}))>(K+\delta)(t_{1}-t_{0}).
\end{equation*}
Let $\tau \in \mathbb{R}$ be defined by 
\begin{equation*}
\tau :=\inf \left\{\,t\in \lbrack t_{0},t_{1}]:~\mathrm{ex\,}(S(t_{0}),S(t))>(K+\delta )(t-t_{0})\,\right\}.
\end{equation*}
By Proposition~\ref{oriented criterion} and the definition of the graphical
modulus $\mathrm{calm}^{\uparrow }$, for each $x\in S(t_{0})$, there is $\varepsilon _{x}>0$ and $\delta_{x}>0$ such that 
\begin{equation*}
S(t_{0})\cap B(x,\delta _{x})\subset S(t)+(K+\frac{\delta}{2})|t-t_{0}|B,\quad \text{for all }t\in \lbrack t_{0},t_{0}+\varepsilon _{x}).
\end{equation*}
Let $\tilde{\varepsilon}_{x}>0$ be the supremum of all $\varepsilon >0$ such
that: 
\begin{equation*}
x\in S(t)+(K+\frac{\delta }{2})|t-t_{0}|B,~\text{for all }t\in \lbrack
t_{0},t_{0}+\varepsilon ).
\end{equation*}
If $\tau =t_{0}$, then there exists a sequence $\{x_{k}\}_{k}\subset S(\tau )
$ such that $\tilde{\varepsilon}_{x_{k}}<1/k$, for all $k\geq 1$. Since $S(\tau )$ is
compact, the sequence $\{x_{k}\}_{k}$ has some cluster point $\overline{x}\in
S(\tau )$. By Proposition~\ref{oriented criterion}, there exist $\varepsilon
_{\overline{x}}>0$ and $\delta _{\overline{x}}>0$ such that 
\begin{equation*}
S(t_{0})\cap B(\overline{x},\delta _{\overline{x}})\subset S(t)+(K+\frac{\delta }{2})|t-t_{0}|B,~\text{for all }t\in \lbrack t_{0},t_{0}+\varepsilon_{\overline{x}}).
\end{equation*}
which contradicts the maximality of $\tilde{\varepsilon}_{x_{k}}$, for $k$
large enough. This establishes that $t_{0}<\tau$. Proceeding in the same way, we can actually show that $\tau \geq t_{1}$. Indeed, assuming $\tau  < t_1$, and using the same argument as above (with $t_{0}$ in the place of $\tau$) together with the triangle inequality we obtain a contradiction in a similar way. Therefore, for every 
$\delta>0$ we have:
\begin{equation*}
\mathrm{ex}(S(t_{0}),S(t_{1}))\leq (K+\delta)(t_{1}-t_{0}),
\end{equation*}
which finishes the first assertion of the lemma. \smallskip \newline
For the second part, we follow the same procedure to estimate the reverse excess 
$\mathrm{ex\,}(S(t_{1}),S(t_{0}))$, and conclude thanks to the fact that $
\mathrm{dist}(S(t_{0}),S(t_{1}))=\max \{\mathrm{ex}(S(t_{0}),S(t_{1})),\mathrm{ex}(S(t_{1}),S(t_{0})) \}$.  The details are left to the
reader.\hfill $\Box $

\bigskip

Now, we proceed with the proof of our second main result.

\subsection{Proof of Theorem B. }

\label{ss: proof of B}

We recall from Section~\ref{ss: proof of A} the definition of $\mathcal{T}$
and fix $a\in \mathcal{T}$. We prove $(a)\Rightarrow (e)\Rightarrow
(f)\Rightarrow (d).\medskip $\newline
$\bm{a)\to e)}:$ Choose $b>a$ such that the statements $(a)$--$(d)$ of
Theorem~A hold true, and $\varphi ^{\uparrow }(t)<+\infty $ for all $t\in
(a,b)$ (\emph{c.f.} Assumption (A2)). We set 
\begin{equation*}
\sigma (t)=\int_{a}^{t}\varphi ^{\uparrow }(s)ds,\quad t\in (a,b).
\end{equation*}
By $(d)$ the above integral is well-defined and $\sigma $ is continuous with 
$\sigma (a)=0.$ Let $\{(t_{i},x_{i})\}_{i\geq 0}\subset \mathcal{S}$ be any
catching-up sequence for $S$ with  $I:=[t_{0},t_{k}]\subset (a,b).$ We shall
prove that \eqref{ortega} holds for every $k\geq 1$. By Proposition~\ref{Prop_seb}, $S$ is continuous on the interval $[t_{0},t_{k}]$ and by Lemma~\ref{continuity of the talweg}, $\varphi ^{\uparrow }$ is continuous (and
finite), hence Riemann integrable there. Let $\{s_{j}^{i}\}_{j=0}^{k_{i}}$ be a
partition of the interval $[t_{i},t_{i+1}],$ $i\in \{0,\ldots ,k-1\},$ with
width
\begin{equation*}
\underset{j\in \{0,\ldots ,k_{i}-1\}}\max \,|s_{j+1}^{i}-s_{j}^{i}|\, <\frac{1}{N},\quad \text{for all }i\in
\{0,\ldots ,k-1\}.
\end{equation*}
Notice that for every $i\in \{0,\ldots ,k-1\},$ we have $s_{0}^{i}=t_{i}$
and $s_{k_{i}}^{i}=t_{i+1}.$ We set 
\begin{equation*}
z_{0}^{i}:=x_{i}\in S(t_{i})\qquad \text{and for each \ }j\in \{0,\ldots
,k_{i}-1\}\;\text{we pick }{z_{j+1}^{i}\in \mathrm{Proj}_{S(s_{j+1})}(z_{j}^{i}).}
\end{equation*}
Then using triangle inequality and the fact that 
\begin{equation*}
 \|x_{i+1} -x_i\|= d(x_i,\underbrace{S(t_{i+1})}_{=S(s^i_{k_i})}))\leq \| z^i_{k_i} -z^i_0 \| .
\end{equation*}
we deduce from Lemma \ref{distance-projection} that:
\begin{equation*}
\Vert x_{i+1}-x_{i}\Vert \leq \sum_{j=0}^{k_{i}-1}\Vert
z_{j+1}^{i}-z_{j}^{i}\Vert \leq \sum_{j=0}^{k_{i}-1}\left( \underset{t\in
\lbrack s_{j}^{i},s_{j+1}^{i}]}{\sup }\varphi ^{\uparrow }(t)\right) \left(
s_{j+1}-s_{j}\right) .
\end{equation*}
Taking the limit as $N\rightarrow \infty $ we obtain that
\begin{equation*}
\Vert x_{i+1}-x_{i}\Vert \leq \int_{t_{i}}^{t_{i+1}}\varphi ^{\uparrow }(s)ds
\end{equation*}
and consequently,
\begin{equation*}
\sum_{i=0}^{k-1}\Vert x_{i+1}-x_{i}\Vert \leq \int_{t_{0}}^{t_{k}}\varphi
^{\uparrow }(t)dt=\sigma (t_{k})-\sigma (t_{0}).
\end{equation*}
\newline
$\bm{e)\to f)}:$ It follows directly by taking $M=\sigma (b)$.\medskip 
\newline
$\bm{f) \to d)}:$ Let $b>a$ and $M>0$ be given by statement $(f)$. By
(A2)--(A3), shrinking $b$ if necessary, we may assume that $\varphi
^{\uparrow }(t)<\infty ,$ for all $t\in (a,b)$ and $\partial S$ is
continuous on $(a,b)$. Notice that for any compact interval $[c,d]\subset
(a,b),$ the function $\varphi ^{\uparrow }$ is continuous and finite on $%
[c,d],$ therefore Riemann integrable. We shall prove that its integral is
bounded by $M$ (independently of the values of $c$ and $d$).\smallskip 
\newline
To this end, let $t_{0}\in \lbrack c,d]$ and $N\in \mathbb{N}$. By
compactness, there exists $x\in S(t_{0})$ such that $\Vert D^{\ast
}S(t_{0},x)|^{+}=\varphi ^{\uparrow }(t_{0})$. If $\varphi ^{\uparrow
}(t_{0})<\frac{1}{N}$, we set $t_{1}:=\min \{t_{0}+\frac{1}{N},d\}$, $x_{0}=x
$ and $y_{0}\in \mathrm{Proj}_{S(t_{1})}(x_{0})$. Observe that 
\begin{equation*}
\Vert x_{0}-y_{0}\Vert \geq 0\geq (t_{1}-t_{0})\left( \varphi ^{\uparrow
}(t_{0})-\frac{1}{N}\right) .
\end{equation*}
If $\varphi ^{\uparrow }(t)\geq \frac{1}{N}$, by Proposition~\ref{oriented
criterion}, since $\mathrm{calm}^{\uparrow }S(t,x)=\varphi ^{\uparrow }(t)$,
there are $x_{0}\in S(t_{0})$ and $t_{1}\in (t_{0},\min \{t_{0}+\frac{1}{N},b\})$ such that any $y_{0}\in \mathrm{Proj}_{S(t_{1})}(x_{0})$ satisfies 
\begin{equation*}
\Vert x_{0}-y_{0}\Vert \geq (t_{1}-t_{0})\left( \varphi ^{\uparrow }(t_{0})-\frac{1}{N}\right) .
\end{equation*}
Using transfinite induction we obtain an increasing net $\{t_{\lambda
}\}_{\lambda \leq \Lambda }\subset \lbrack c,d]$ indexed over a countable
ordinal $\Lambda $, such that $t_{0}=c$, $t_{\Lambda }=d$, $0<t_{\lambda+1}-t_\lambda\leq 1/N$ for all $\lambda<\Lambda$, and for any limit ordinal $\alpha\leq\Lambda$, $t_\alpha=\sup\{t_\lambda:~\lambda<\alpha\}$. Also, we get a net 
$\{(x_{\lambda },y_{\lambda })\}_{\lambda \leq \Lambda }$ such that $x_\lambda\in S(t_\lambda)$, $y_\lambda\in \textup{Proj}_{S(t_{\lambda+1})}(x_\lambda)$ and
\begin{equation*}
\Vert x_{\lambda }-y_{\lambda }\Vert \geq (t_{\lambda +1}-t_{\lambda
})\left( \varphi ^{\uparrow }(t_{\lambda })-\frac{1}{N}\right) ,~\text{for
all }\lambda <\Lambda.
\end{equation*}
For every finite subset $F\subset \Lambda $ we have 
\begin{equation*}
\sum_{\lambda \in F}\Vert x_{\lambda }-y_{\lambda }\Vert \geq \left (\sum_{\lambda
\in F}(t_{\lambda +1}-t_{\lambda })\varphi ^{\uparrow }(t_{\lambda })\right)-\dfrac{d-c}{N}.
\end{equation*}
Since $\{(t_{\lambda},x_\lambda),~(t_{\lambda},y_\lambda):~\lambda\in F\}$ is a subsequence of a piecewise catching-up sequence for $S$, 
taking the supremum over all finite families $F$ of $\Lambda $ we get
\begin{equation*}
M\geq \sum_{\lambda < \Lambda }\Vert x_{\lambda }-y_{\lambda }\Vert \geq \left(
\sum_{\lambda < \Lambda }(t_{\lambda +1}-t_{\lambda })\varphi ^{\uparrow
}(t_{\lambda }) \right)-\dfrac{d-c}{N}.
\end{equation*}
Taking the limit as $N$ goes to infinity we obtain:
\begin{equation*}
M\geq \int_{c}^{d}\varphi ^{\uparrow }(t)dt.
\end{equation*}
Since the above is independent of the interval $[c,d],$ we deduce that $\varphi ^{\uparrow }$ is integrable on $(a,b)$.\hfill $\Box $

\bigskip

\noindent\textbf{Acknowledgements.} The first author thanks J. Bolte (Toulouse) for useful discussions which originated this work.



\noindent\rule[0pt]{5cm}{2pt}

\vspace{0.5cm}

\noindent Aris DANIILIDIS

\medskip

\noindent DIM--CMM, UMI CNRS 2807\newline Beauchef 851, FCFM, Universidad de
Chile \smallskip

\noindent E-mail: \texttt{arisd@dim.uchile.cl}
\newline\noindent\texttt{http://www.dim.uchile.cl/\symbol{126}arisd}

\medskip

\noindent Research supported by the grants: \newline CMM AFB170001, FONDECYT
1211217 (Chile), ECOS-ANID C18E04 (Chile, France).\\
\vspace{0.6cm}

\noindent Sebastián TAPIA GARCIA

\medskip

\noindent DIM--CMM, UMI CNRS 2807\newline Beauchef 851, FCFM, Universidad de
Chile \newline\smallskip

\noindent IMB, UMR CNRS 5251\newline
Cours de la Lib\'eration 351, Talence, Université de Bordeaux.

\noindent E-mail:  \texttt{stapia@dim.uchile.cl}

\medskip

\noindent Research supported by the grants: \newline CMM AFB170001, FONDECYT
1211217 (Chile), ECOS-ANID C18E04 (Chile, France), ANID-PFCHA/Doctorado Nacional/2018-21181905.
\end{document}